\documentclass[11pt]{article}
\usepackage[auth-sc]{authblk}
\usepackage{amssymb}
\usepackage{graphicx}
\usepackage{latexsym}
\usepackage{amsmath,amscd}
\usepackage{bm}
\def\part#1{\frac{\partial\phantom{q}}{\partial#1}}

\newenvironment{rmk}{\begin{trivlist}\item[]{\bf Remark:} }
{\end{trivlist}}
\newenvironment{rmks}{\begin{trivlist}\item[]{\bf Remarks:} }
{\end{trivlist}}

\newenvironment{prf}{\begin{trivlist}\item[]{\bf Proof:} }
{\hfill $\Box$ \end{trivlist}}

\newtheorem{thm}{Theorem}

\newtheorem{prp}[thm]{Proposition}

\newcommand{\lie}[1]{\mathfrak{#1}}
\def\End{\mathop{\rm End}\nolimits}

\def\deg{\mathop{\rm deg}\nolimits}

\def\rk{\mathop{\rm rk}\nolimits}
\def\tr{\mathop{\rm tr}\nolimits}

\newcommand{\R}{\mathbf{R}}
\newcommand{\C}{\mathbf{C}}
\newcommand{\K}{\mathbf{H}}
\newcommand{\Z}{\mathbf{Z}}

\newcommand{\CP}{{\mathbf P}}

\title{Nonabelianization of Higgs bundles}
 \author[1]{Nigel Hitchin}
 \author[2]{ Laura P. Schaposnik}
 \affil [1] {Mathematical Institute, 24-29 St Giles, Oxford OX1 3LB, UK}
 \affil [2] {Mathematisches  Institut, Ruprecht-Karls-Universit\"at, 69120 Heidelberg, Germany}

\begin{document}
 \maketitle
 \thispagestyle{empty}

\let\oldthefootnote\thefootnote
\renewcommand{\thefootnote}{\fnsymbol{footnote}}
\footnotetext{hitchin@maths.ox.ac.uk \hskip 1cm lauraschaposnik@gmail.com}
\let\thefootnote\oldthefootnote

\section{Introduction}
The moduli space of Higgs bundles over a compact Riemann surface $\Sigma$ of genus $g>0$ for a complex group $G^c$  has the well-known structure of a completely integrable Hamiltonian system: a proper map to an affine space, whose generic fibre is an abelian variety. For the general linear group  such a Higgs bundle consists of   a  vector bundle $V$ together with a section $\Phi$ of $\End V\otimes K$. The coefficients of the characteristic polynomial of $\Phi$ define the map to the base and the corresponding  fibre is the Jacobian of an algebraic curve $S$ defined by the equation $\det (xI-\Phi)=0$. This is a covering  $\pi:S\rightarrow \Sigma$ on which $\Phi$ has a single-valued eigenvalue $x$, a section of $\pi^*K$. Conversely,  given a line bundle $L$ on $S$ we obtain $V$ as the direct image sheaf of $L$ and $\Phi$ as the direct image of $x:L\rightarrow L\otimes \pi^*K$. This  {\it abelianization} process has been useful in attacking various problems relating to bundles on curves.

It is clear, however, that we could replace $L$ by a rank $r$ bundle $E$ on $S$ and obtain a Higgs bundle on $\Sigma$ by the same construction. In this case each generic eigenspace of $\Phi$ is $r$-dimensional and $\det (xI-\Phi)=p(x)^r$ for some polynomial $p(x)$. Since the map to the base is surjective such  cases certainly occur. What we show here is that they occur, with $r=2$, very naturally when considering the Higgs bundles which correspond (by solving the gauge-theoretic Higgs bundle equations) to flat connections on $\Sigma$ with holonomy in  the real Lie groups $G^r=SL(m,\K), SO(2m,\K)$ and $Sp(m,m)$, real forms of $G^c= SL(2m,\C), SO(4m,\C)$ and $Sp(4m,\C)$ respectively. The first two groups are often denoted by $SU^*(2m)$ and $SO^*(2m)$ but the phenomenon we are describing clearly reflects the noncommutativity of the quaternions which justifies the former  notation.

 The fibres of the integrable system, even those over non-regular values,  are always compact. We find here that for $SL(m,\K)$ the fibre consists of the moduli space of semi-stable rank $2$ bundles  with fixed determinant on the spectral curve $S$. For  $SO(2m,\K)$ the fibre has many components each of which is a moduli space of semi-stable rank $2$ bundles on a quotient $\bar S$ of the spectral curve and for $Sp(m,m)$ it is  a $\Z_2$-quotient of a moduli space of semi-stable rank $2$ {\it  parabolic} bundles on $\bar S$.

Part of this work is contained in the second author's DPhil thesis  \cite{LS}. She acknowledges the support of Oxford University, New College and QGM  Aarhus.  Both authors wish to thank ICMAT, Madrid  for support during an activity of the ``Nigel Hitchin Laboratory"   in March 2013 where this article originated. Work of Ana Pe\'on, to appear in a forthcoming thesis,  was also presented there, and gives a characterization using cameral covers of the real forms which exhibit this nonabelianization phenomenon.  

\section{Higgs bundles for real forms}
Under suitable stability conditions a Higgs bundle defines a solution of  equations for a $G$-connection $A$, where $G$ is the maximal compact subgroup of $G^c$ \cite{Hit1},\cite{Sim}. For $G=U(n)$ these are $F_A+[\Phi,\Phi^*]=0$ and then the connection $\nabla_A+\Phi+\Phi^*$ is flat, with holonomy in $GL(n,\C)$. To get a flat connection with holonomy in $GL(n,\R)$ we take $A$ to be an $SO(n)$-connection and $\Phi=\Phi^T$ using the transpose defined by the orthogonal structure on $V$. 

For a  general real form $G^r$ of $G^c$ we take a $U$-connection where $U$ is the maximal compact subgroup of $G^r$ and in the decomposition $\lie{g}=\lie{u}\oplus \lie{m}$ take $\Phi\in H^0(\Sigma,\lie{m}\otimes K)$. Much work on the study of connected components of moduli spaces of flat $G^r$-connections using this approach has been carried out by Bradlow, Garcia-Prada et al. (for example in \cite{BGP},\cite{GGM}). For the groups in question we have the following descriptions of the corresponding Higgs bundles:

\begin{itemize}
\item
The group $SL(m,\K)$  is the subgroup of $SL(2m,\C)$ which commutes with an antilinear automorphism $J$ of $\C^{2m}$ such that $J^2=-1$. Its maximal compact subgroup is the  quaternionic unitary group $Sp(m)$. Since $Sp(m)^c=Sp(2m,\C)$, the corresponding Higgs bundle  consists of a rank $2m$ symplectic vector bundle $(V,\omega)$ and the  Higgs field satisfies $\Phi=\Phi^T$ for the symplectic transpose. Using $\omega$  to identify $V$ and $V^*$, this means that $\Phi=\omega^{-1}\phi$ for $\phi\in H^0(\Sigma,\Lambda^2V\otimes K)$. 
\item
The group $SO(2m,\K)$ is the subgroup of $GL(4m,\C)$ which preserves a complex inner product $(u,v)$ and commutes with an antilinear automorphism $J$ as above for which $(u,Jv)$ is a hermitian form. If $(u,Ju)>0$ then $(Ju,J^2u)=-(Ju,u)<0$ and so the form has hermitian signature $(2m,2m)$.  The maximal compact subgroup is $U(2m)$ and the Higgs bundle has the form $V=W\oplus W^*$ for a rank $2m$ vector bundle $W$. The inner product is defined by the natural pairing between $W$ and $W^*$. The Higgs field is of the form
$\Phi(w,\xi)=(\beta(\xi),\gamma(w))$ where $\beta:W^*\rightarrow W\otimes K$ and $\gamma:W\rightarrow W^*\otimes K$ are skew-symmetric. 
\item
The group  $Sp(m,m)$ is the subgroup of $GL(4m,\C)$ which preserves a complex symplectic form  $\omega$ and commutes with an antilinear automorphism $J$ for which   $\omega( u,Jv)$ is a Hermitian form of signature $(2m,2m)$. It is the group of quaternionic matrices which are unitary with respect to an indefinite form. The maximal compact subgroup is $Sp(m)\times Sp(m)$ and the Higgs bundle is of the form $W_1\oplus W_2$ for symplectic rank $2m$ vector bundles $(W_1,\omega_1),(W_2,\omega_2)$. The Higgs field is of the form
$\Phi(v,w)=(\beta(w),\gamma(v))$ where $\beta:W_2\rightarrow W_1\otimes K$ and $\gamma:W_1\rightarrow W_2\otimes K$ with $\beta=-\gamma^T$, using the symplectic transpose. 
\end{itemize}
Since $SO(2m,\K)$ and $Sp(m,m)$ are subgroups of $SL(2m,\K)$ we begin dealing with the first case.

\section{Spectral data for $SL(m,\K)$}\label{SL}

As noted above, in this case we have a rank $2m$ symplectic vector bundle $(V,\omega)$ and a Higgs field which is symmetric with respect to $\omega$. If $A$ is a symmetric endomorphism of a symplectic vector space $U$ of dimension $2m$, and $\alpha$ the corresponding element of $\Lambda^2U^*$, then $xI-A$ is singular if and only if the exterior product $(x\omega-\alpha)^m=0$. The Pfaffian polynomial of $xI-A$ is defined by $p(x)\omega^m=(x\omega-\alpha)^m$, and then $p(x)^2=\det(xI-A)$. If $p(x)$ has distinct roots then $U$ decomposes into a sum of $m$ two-dimensional symplectic eigenspaces of $A$ and in this case, and by continuity in all cases, $A$ satisfies the matrix equation $p(A)=0$.  If $A$ is in the Lie algebra of $SL(2m,\C)$ then in addition $\tr A=0$. 

Replacing $A$ by $\Phi$ we have a polynomial  $p(x)=x^m+a_2x^{m-2}+\dots+a_m$ where the coefficients $a_i\in H^0(\Sigma, K^i)$. 
As an $SL(2m,\C)$ Higgs bundle, the usual spectral curve is defined by the vanishing of the characteristic polynomial   
 so the coefficients of $p(x)^2$, which lie in $\bigoplus_{i=2}^{2m}H^0(\Sigma,K^i)$, define a point in the base of the fibration. In our case, by a slight abuse of notation,  we shall call the curve $S$ defined by $p(x)=0$ the spectral curve. Bertini's theorem assures us that for generic $a_i$ the curve is nonsingular. It is a ramified $m$-fold cover of $\Sigma$. More precisely, we may interpret the equation $p(x)=0$ as the vanishing of a section of $\pi^*K^m$ over the total space of the canonical bundle $\pi:K\rightarrow \Sigma$, where $x$ is the tautological section of $\pi^*K$. The cotangent bundle of $\Sigma$ is a symplectic manifold and hence has trivial canonical bundle, so $K_S\otimes \pi^*K^{-m}$ is trivial and  $K_S= \pi^*K^{m}$. Taking degrees of both sides this says that the genus  of $S$ is given  by $g_S=m^2(g-1)+1.$

 On the spectral curve $S$, $x$ is a well-defined eigenvalue of  $\Phi$, and the cokernel of $xI-\Phi$ is a rank two holomorphic vector bundle $E$. It then follows, as in \cite{BNR} (and using $p(\Phi)=0$ instead of the Cayley-Hamilton theorem), that we can identify $V$ with the direct image $\pi_*E$ and the Higgs field $\Phi$  as the direct image of $x:E\rightarrow E\otimes \pi^*K$ (recall that the direct image sheaf is defined for each open set $U\subset \Sigma$ by $H^0(U,\pi_*E)=H^0(\pi^{-1}(U),E)$).

 If we now start with any rank $2$ bundle $E$ on $S$ we can obtain by the same construction a $GL(2m,\C)$ Higgs bundle, but we need to determine the conditions on $E$ for this to be the data for the group $SL(m,\K)$:
 
  \begin{prp} \label{prop1}Let $p=x^m+a_2x^{m-2}+\dots+a_m$ be a section of the line bundle  $\pi^*K^m$ on the cotangent bundle of $\Sigma$ whose divisor is a smooth curve $S$,  and let $E$ be a rank $2$ vector bundle on $S$. Then the  direct image of $x:E\rightarrow E\otimes \pi^*K$    defines a semi-stable Higgs bundle on $\Sigma$ for the group $SL(m,\K)$ if and only if  
  \begin{itemize}
  \item 
  $\Lambda^2E\cong \pi^*K^{m-1}$ 
  \item
  $E$ is semi-stable.
  \end{itemize}
    \end{prp}
 \begin{prf} First we define a nondegenerate skew form on $\pi_*E$. For this, the relative duality theorem gives 
 $$(\pi_*E)^*\cong \pi_*(E^*\otimes K_S)\otimes K^{-1}$$
so if $V=\pi_*E$, to achieve $V\cong V^*$ we want $E^*\otimes \pi^*K^{m-1}\cong E$. But  we are given $\Lambda^2 E\cong\pi^*K^{m-1}$ so this is satisfied. We should  also check that the duality is provided by a skew form, and this requires a  concrete expression of relative duality. At a regular value $a\in \Sigma$ of $\pi$, 
$$V=\bigoplus_{y\in \pi^{-1}(a)} E_y$$ 
and taking $s_y\in E_y,s'\in E^*_y\otimes K_S\otimes \pi^{*}K^{-1}$ we form 
$$\sum_{y\in \pi^{-1}(a)}\frac{\langle s,s'\rangle}{d\pi_y}$$
where $d\pi:K^{-1}_S\rightarrow \pi^*K^{-1}$ is the derivative, considered in $K_S\otimes \pi^*K^{-1}$. This extends for the direct image over branch points. But the isomorphism  $E^*\otimes \pi^*K^{m-1}\cong E$ is skew-symmetric, showing that $V\cong V^*$ is also skew. Moreover multiplication by $x$ satisfies $\langle x s,s'\rangle=\langle s,x s'\rangle$ which is symmetric and defines $\Phi$ as a Higgs field satisfying $\Phi=\Phi^T$.
\vskip .25cm
The semi-stability condition for Higgs bundles \cite{Hit1} is that, for each $\Phi$-invariant subbundle $W\subset V$ , $\deg W/\rk W\le \deg V/\rk V$. But $V$ is symplectic so  $\deg V=0$, and then we require $\deg W\le 0$. Now  $\Phi\vert_W$ has a characteristic polynomial which, if $W\ne V$, divides that of $\Phi$. But the characteristic polynomial of $\Phi$  is $p(x)^2$ and $S$ is smooth, and in particular irreducible. So the characteristic polynomial for $\Phi\vert_W$ must be $p(x)$. Then $(W,\Phi)$ is a rank $m$ Higgs bundle and by \cite{BNR} is the direct image of a line bundle $L\subset E$. 

From Grothendieck-Riemann-Roch we have   $(1-g)m+\deg W=(1-g_S)+\deg L$
and so $\deg W=(1-g)(m^2-m)+\deg L$. Hence the Higgs semi-stability condition is $\deg L\le m(m-1)(g-1)$. But $\deg E=\deg \pi^*K^{m-1}=m(m-1)(2g-2)$ and so this is equivalent to the semi-stability condition $\deg L\le \deg E/2$ for the rank $2$ bundle $E$.
 \end{prf}
 
 \begin{rmks}
 
 \noindent 1. The degree of $E$ on $S$ is $2m(m-1)(g-1)$ which is even and so the moduli space of semi-stable bundles is singular, the singular locus represented by decomposable bundles $E=L\oplus (L^*\otimes \pi^*K^{m-1})$. Then relative duality gives $\pi_*E=W\oplus W^*$ and the Higgs field is $\Phi=(\phi,\phi^T)$ for a $GL(m,\C)$ Higgs bundle $(W,\phi)$.

 \noindent 2.  One may check dimensions of the moduli space here: considered as the moduli space of a real form its real dimension is $(2g-2)\dim G^r=2(4m^2-1)(g-1)$. The complex dimension of the space of polynomials $p(x)$ is $3(g-1)+5(g-1)+\dots +(2m-1)(g-1)=(m^2-1)(g-1)$, and for the moduli space of stable bundles on $S$ with fixed determinant it is $3(g_S-1)=3m^2(g-1)$ giving in total $(4m^2-1)(g-1)$.

\end{rmks}

 \section{Spectral data for $SO(2m,\K)$}\label{SO}
 
 The Higgs bundle here is $V=W\oplus W^*$ where $W$ has rank $2m$ and the Higgs field is of the form 
 \begin{equation}
 \Phi=\begin{pmatrix}0 & \beta\\
 \gamma & 0\end{pmatrix}
 \label{off}
 \end{equation}
 where  $\beta:W^*\rightarrow W\otimes K$ and $\gamma: W\rightarrow W^*\otimes K$ are both skew-symmetric. The inclusion $SO(2m,\K)\subset SL(2m,\K)$ means we should also consider this as a special case of the previous section. To do this, define a symplectic form on $V$ by
 $\omega((w_1,\xi_1),(w_2,\xi_2))=\xi_2(w_1)-\xi_1(w_2)$. Then
 $$\omega(\Phi(w_1,\xi_1),(w_2,\xi_2))=\xi_2(\beta\xi_1)-\gamma w_1(w_2)=-\xi_1(\beta \xi_2)+\gamma w_2 (w_1)=\omega((w_1,\xi_1),\Phi(w_2,\xi_2))$$
 and so $\Phi=\Phi^T$.
 
 It follows that  we can use the result of the previous section to deduce that $\det (xI-\Phi)=p(x)^2$ for some polynomial of degree $2m$ and, assuming it defines a smooth curve $S$, $\Phi$ has generically two-dimensional eigenspaces. Globally, as before  the bundle $V$ can be written $\pi_*E$ for a rank $2$ bundle on the spectral curve $S$, which has  genus $g_S=4m^2(g-1)+1$. 
 
 Suppose that $(w,\xi)\in W\oplus W^*$ is an eigenvector of $\Phi$ with eigenvalue $\lambda$. Then $\beta(\xi)=\lambda w$ and $\gamma(w)=\lambda \xi$. Hence  
 $$\Phi(w,-\xi)=(-\lambda w,\lambda \xi)=-\lambda(w,-\xi).$$
 Thus for each two-dimensional generic eigenspace with eigenvalue $\lambda$ there exists another with eigenvalue $-\lambda$. In particular this means that 
 $p(x)=x^{2m}+a_1x^{2m-2}+\dots +a_m$ where $a_i\in H^0(\Sigma, K^{2i})$ and the curve $S$ has an involution $\sigma(x)=-x$. We now need to determine the properties of  the rank $2$ bundle $E$:
 
  \begin{prp} \label{prop2} Let $p=x^{2m}+a_2x^{2m-2}+\dots+a_m$ be a section of the line bundle  $\pi^*K^{2m}$ on the cotangent bundle of $\Sigma$ whose divisor is a smooth curve $S$ and let  $\sigma$ be the involution  $\sigma(x)=-x$. If $E$ is a rank $2$ vector bundle on $S$, then the  direct image of $x:E\rightarrow E\otimes \pi^*K$    defines a semi-stable Higgs bundle on $\Sigma$ for the group $SO(2m,\K)$ iff
   \begin{itemize}
  \item 
  $\Lambda^2E\cong \pi^*K^{m-1}$ 
  \item
  $E$ is semi-stable
  \item
  $\sigma^*E\cong E$ where the induced action on  $\Lambda^2E=\pi^*K^{2m-1}$ is trivial.
  \end{itemize}
    \end{prp}
    
    \begin{rmk} The bundle $\pi^*K^{2m-1}$ is pulled back from $\Sigma$ and this is why it makes sense to speak of the trivial action of $\sigma$, since $\pi\sigma=\pi$. When $E$ is stable any automorphism is a scalar so the isomorphism gives a lifted action of the involution, well-defined modulo $\pm 1$ and in particular the action on $\Lambda^2E$ is well-defined. In fact in general all we require is that the action on $\Lambda^2E$ at fixed points should be trivial. 
    \end{rmk}
    
    \begin{prf}  
    The first part follows from Proposition \ref{prop1}. The bundle $E$ is defined as the cokernel of $xI-\Phi$ in $V\otimes K$, but using the orthogonal structure on $V$, this means that $E^*\otimes \pi^*K $  is the kernel of $xI+\Phi$, i.e. generically the two-dimensional eigenspace of $\Phi$ with eigenvalue $-x$. We saw above how $(w,\xi)\mapsto (w,-\xi)$ gives an isomorphism between the $\pm x$    eigenspaces and so there is an isomorphism $\sigma^*E\cong E$. 
    
Triviality of the action on $\Lambda^2E=\pi^*K^{2m-1}$ in the statement of the Proposition is a function of the isomorphism, which came to us  as in Section \ref{SL}   from relative duality.  This gave us the following formula for  the symplectic form   over a regular value: 
     \begin{equation}
     \sum_{y\in \pi^{-1}(a)}\frac{ s\wedge s'}{d\pi_y}
     \label{sympW}
     \end{equation}
    where  $d\pi$, the derivative, is a section of $K_S\otimes  \pi^*K^{-1}$ and we use the canonical symplectic form  of the cotangent  bundle of $\Sigma$ to identify this with $\pi^*K^{2m-1}$. Together with an isomorphism $\pi^*K^{2m-1}\cong \Lambda^2 E$,  formula (\ref{sympW}) is a well-defined scalar. Now the  symplectic form on the cotangent bundle is anti-invariant under $\sigma$, scalar multiplication by $-1$ in the fibres, and the pairing symplectic form on $W\oplus W^*$ is anti-invariant under the map $(w,\xi)\mapsto (w,-\xi)$.  It follows that the action on $\Lambda^2E=\pi^*K^{2m-1}$ is trivial.
     
      Now let $\bar S$ be the quotient of $S$ by the involution and $p:S\rightarrow \bar S$ the quotient map. Setting $z=x^2$ embeds $\bar S$ in the total space of $K^2$ with equation $z^{m}+a_1z^{m-1}+\dots +a_m=0$. If $\bar \pi:K^2\rightarrow \Sigma$ is the projection then $z$ is the tautological section of $\bar\pi^*K^2$ and $\pi:S\rightarrow \Sigma$ can be written as $\pi=\bar\pi\circ p$.
    
     For an open set $U\subset \bar S$, $p^{-1}(U)$ is $\sigma$-invariant and then the $\pm 1$ eigenspaces   of the action on $H^0(p^{-1}(U),E)$ decompose   the direct image on $\bar S$ as  $p_*E=E^+\oplus E^-$.  Since $(w,0)\in W$ and $(0,\xi)\in W^*$ are the $\pm 1$ eigenspaces of the involution we see that $\bar\pi_*E^+=W,\bar\pi_*E^-=W^*$.

          \vskip .25cm
     Conversely, suppose we are given $S$ and $E$ as in the statement of the Proposition. Then $\pi_*E=\bar\pi_*(E^+\oplus E^-)=W_1\oplus W_2$ is a symplectic vector bundle from Proposition \ref{prop1}. Moreover since $\sigma(x)=-x$, $x:E\rightarrow E\otimes \pi^*K$ maps local invariant sections of $E$ to anti-invariant ones and so the Higgs field has the off-diagonal shape of (\ref{off}). 
     To obtain $V=W\oplus W^*$ we need to show that $W_1$ and $W_2$ are Lagrangian with respect to the symplectic structure on $V$. Now $W_1$ is defined by the direct image of an invariant section $s$. Equation (\ref{sympW}) evaluates the symplectic pairing of two sections $s,s'$  but if they are both invariant then $s\wedge s'$ is invariant. As in the discussion above, if the action on $\Lambda^2E$ is trivial, then the denominator in   (\ref{sympW}) is anti-invariant and  so the terms over $y$ and $\sigma(y)$ cancel and so the symplectic form on $V$ vanishes on $W_1$. A similar argument holds for $W_2$.   They are therefore transverse Lagrangian subbundles and setting $W_1=W$, $W_2\cong W^*$. 
    
    As a special case of the previous section, $\Phi=\Phi^T$ with respect to the symplectic form but given its shape (\ref{off}) this means that the terms $\beta$ and $\gamma$ are skew-symmetric. 
    \end{prf}
    
Proposition \ref{prop2} tells us that the fibre of the integrable system is defined by the fixed points of an involution induced by $\sigma$ on the moduli space of rank $2$ semi-stable bundles on $S$. There are several components, however. This is clear from the flat connection point of view: the maximal compact subgroup of  $SO(2m,\K)$ is $U(2m)$ and so any flat $SO(2m,\K)$ bundle can be topologically reduced to $U(2m)$ where it has a Chern class. In the Higgs bundle description this is the degree of the vector bundle $W$. As in the case of $U(m,m)$ dealt with in \cite{LS} we can determine this invariant by considering the action at the fixed points of $\sigma$ on $S$. 
 
 At a fixed point $a$ of $\sigma$ there is a linear  action of $\sigma$ on the fibre $E_a$. Since the action on $\Lambda^2E_a$ is trivial this is scalar multiplication $\pm 1$ and we can assign to each fixed point this number. 
 
 \begin{prp} Suppose the action is $+1$ at $M$ fixed points, then $\deg W=2M-4m(g-1)$.
  \end{prp}
 
 \begin{prf} The fixed point set of $\sigma$ is the intersection of the zero section of $K$ with $S$. Setting $x=0$ in the equation $x^{2m}+a_1x^{2m-2}+\dots +a_m=0$, these points are the images of the $4m(g-1)$ zeros of $a_m\in H^0(\Sigma,K^{2m})$ under the zero section. The action is $+1$ at $M$ of these points. 
 
  Choose a line bundle $L$ on $\Sigma$ of large enough degree such that $H^1(\Sigma, V\otimes L)=0$.   
   By definition of $E^+,E^-$ we have 
 $\dim H^0(S,E\otimes \pi^*L)^{\pm}=\dim H^0(\bar S,E^{\pm}\otimes \bar\pi^*L)$ where the superscript denotes the $\pm$ eigenspace under the action of $\sigma$.  Since $V=\pi_*E$ and  $H^1(\Sigma, V\otimes L)=0$ the higher cohomology groups vanish and 
 applying the holomorphic Lefschetz formula  we obtain
 $$\dim H^0(\bar S,E^{+}\otimes \bar\pi^*L)-\dim H^0(\bar S,E^{-}\otimes \bar\pi^*L)=2(M-(4m(g-1)-M))$$
 and Riemann-Roch gives 
  $$\dim H^0(\bar S,E^{+}\otimes \bar\pi^*L)+\dim H^0(\bar S,E^{-}\otimes \bar\pi^*L)=\dim H^0(\Sigma, V\otimes L)=4m(1-g+\deg L).$$ since $V$ is symplectic and $\deg V=0$.

 Now $W=\bar\pi_*E^+$, so $\dim H^0(\bar S, E^{+}\otimes \bar\pi^*L) = \dim H^0(\Sigma, W\otimes L)=2m((1-g)+\deg L)+\deg W$ by Riemann-Roch and from these three equations we obtain $\deg W=2M-4m(g-1)$.  
 \end{prf} 
 
 \begin{rmks} 
 
 \noindent 1. Since $M\le 4m(g-1)$ we have $\vert \deg W\vert\le 4m(g-1)$ which is the Milnor-Wood inequality for the group $SO(2m,\K)$. 
 
 \noindent 2. In the maximal case $\deg W=4m(g-1)$ all fixed points have action $+1$ and then the bundle $E$ is pulled back from the curve $\bar S$. In this case $\gamma:W\rightarrow W^*\otimes K$ is a homomorphism of bundles of the same degree and so is either everywhere singular or an isomorphism. But $S$ is smooth so $a_m$ is not identically zero and hence $W\cong W^*\otimes K$, or setting $U=W\otimes K^{-1/2}$, $\Psi=\beta\gamma$ we have a Higgs bundle of the same type as an $SL(m,\K)$ bundle but with a $K^2$-twisted Higgs field $\Psi$. Moreover the spectral data of the previous section holds if one takes the rank $2$ bundle $E\otimes \bar\pi ^*K^{-1/2}$ on $\bar S$. This is a case of the Cayley correspondence of \cite{GGM}.

 \noindent 3. For each choice of $M$ fixed points ${a_1,\dots,a_M}\in S$, using the Narasimhan-Seshadri theorem we can interpret the moduli space of invariant semi-stable rank $2$ bundles on $S$ as the moduli space of representations of the fundamental group of $\bar S\backslash\{{p(a_1),\dots,p(a_M)}\}$ with holonomy $-1$ around the marked points. If $M$ is odd this is the (smooth and connected) moduli space of stable rank $2$ bundles on $\bar S$ of odd degree and fixed determinant and if $M$ is even it is the singular moduli space of bundles of even degree.  Given this we can check dimensions as before. The curve $S$ gives $3(g-1)+7(g-1)+\dots + (4m-1)(g-1)=m(2m+1)(g-1)$ complex parameters, and the moduli space of bundles on $\bar S$ gives $3(g_{\bar S}-1)=3m(2m-1)(g-1)$ parameters.  In total this makes 
$3m(2m-1)(g-1)+m(2m+1)(g-1)=m(8m-2)(g-1)=\dim SO(4m,\C)(g-1).$

 \noindent 4. In the stable case the only other action on $E$ is to multiply the given action by $-1$ which changes $M$ to $4m(g-1)-M$ and interchanges the roles of $W$ and $W^*$. Thus there are $2^{4m(g-1)-1}$ components in the fibre.

 \end{rmks} 
 
  \section{Spectral data for $Sp(m,m)$}
 For this group, the Higgs bundle $V=W_1\oplus W_2$ for symplectic rank $2m$ vector bundles $(W_1,\omega_1),(W_2,\omega_2)$. The Higgs field is \begin{equation}
 \Phi=\begin{pmatrix}0 & \beta\\
 -\beta^T & 0\end{pmatrix}
 \label{off2}
 \end{equation}
 where $\beta^T:W_1\rightarrow W_2\otimes K$ is the symplectic adjoint. Since $Sp(m,m)\subset SL(2m,\K)$ we can  apply the results of Section \ref{SL} but we need a symplectic structure on $V$ for which $\Phi=\Phi^T$. Define $\omega=(\omega_1,-\omega_2)$. Then
 $$\omega(\Phi(u_1,u_2),(w_1,w_2))=\omega_1(\beta u_2,w_1)+\omega_2(\beta^T u_1,w_2)= \omega_1(u_2,\beta^Tw_1)+\omega_2(u_1,\beta w_2)$$
and this is equal to
$$-\omega_1(\beta^Tw_1,u_2)-\omega_2(\beta w_2, u_1)=-\omega(\Phi(w_1,w_2),(u_1,u_2))=\omega((u_1,u_2),\Phi(w_1,w_2)).$$

 Here the only difference with the previous case is the action of the involution on $E$:
 
  \begin{prp} \label{prop3}   Let $p=x^{2m}+a_2x^{2m-2}+\dots+a_m$ be a section of the line bundle  $\pi^*K^{2m}$ on the cotangent bundle of $\Sigma$ whose divisor is a smooth curve $S$ and let  $\sigma$ be the involution  $\sigma(x)=-x$. If $E$ is a rank $2$ vector bundle on $S$, then the  direct image of $x:E\rightarrow E\otimes \pi^*K$    defines a semi-stable Higgs bundle on $\Sigma$ for the group $Sp(m,m)$ iff   \begin{itemize}
  \item 
  $\Lambda^2E\cong \pi^*K^{m-1}$ 
  \item
  $E$ is semi-stable
  \item
  $\sigma^*E\cong E$ where the induced action on  $\Lambda^2E=\pi^*K^{2m-1}$ is $-1$.
  \end{itemize}
    \end{prp}
 \begin{prf} The proof proceeds exactly as in Proposition \ref{prop2} until the point where we prove that $W_1$ and $W_2$ are Lagrangian. With the opposite action on $\Lambda^2E$ we deduce instead that $W_1$ and $W_2$ are symplectically orthogonal and hence $V$ is the symplectic sum of  
 $W_1$ and $W_2$.
 
 A slightly different and more detailed approach may be found in \cite{LS}. 
 \end{prf}
 
\begin{rmks}

\noindent 1.  Given the action of $-1$ on $\Lambda^2E$, at a fixed point we have distinct $+1$ and $-1$ eigenspaces. Following \cite{AG} this defines a rank $2$ bundle on the curve $\bar S$ with a parabolic structure at the fixed points defined by the flag given by the $-1$ eigenspace and the parabolic weight $1/2$.  As in the previous case the choice of action corresponds to an ordering of $W_1$ and $W_2$ so a point in the moduli space for $Sp(m,m)$ determines a point in the quotient of the moduli space of parabolic structures by interchanging the roles of the $+1$ and $-1$ eigenspaces. 
 
\noindent 2. Note that  the two groups $SO(m,\K)$ and $Sp(m,m)$ correspond to the two equivariant structures on the line bundle $\Lambda^2E\cong \pi^*K^{2m-1}$ and, as in \cite{AG},  account for all the fixed points in the moduli space of rank $2$ bundles over $S$.

\noindent 3. We can use the parabolic aspect as a check on the dimension: the parameters for the spectral curve and bundle gives $m(8m-2)(g-1)$ as in the previous case but there is a contribution of $1=\dim \CP^1$ for each of the  $4m(g-1)$ parabolic points giving in total $m(8m+2)(g-1)=\dim Sp(4m,\C)(g-1).$ 
 \end{rmks}
 
 \section{Comments}
 
 \noindent 1. The representation of the moduli space of flat connections as the Higgs bundle moduli space, and in particular the integrable system,  depends on the choice of a complex structure on the underlying real surface $\Sigma$. Properties of the Higgs bundle can change significantly for the same representation of $\pi_1(\Sigma)$. As an example, the uniformizing representation in $PSL(2,\R)$ of a Riemann surface has a nilpotent Higgs field for the natural complex structure but the same representation has a non-singular spectral curve $x^2-a=0$  when we change the complex structure of $\Sigma$ and hence its Higgs bundle moduli space.  It is natural to ask which representations have smooth, or irreducible,  spectral curves in {\it some} complex structure. The examples given here have reducible spectral curves in {\it any} complex structure. 
 
 \noindent 2. The classical abelianization picture of a spectral curve with a line bundle over it is, in the physicists' terminology, a D-brane: a Lagrangian submanifold of the cotangent bundle of $\Sigma$ together with a flat line bundle over it. However, when two such D-branes coalesce one expects to find a flat higher rank bundle over the resulting curve. Given the stability property of the rank $2$ bundle here, it follows from the Narasimhan-Seshadri theorem that this is precisely what we have. A sequence of points in the $GL(2m,\C)$ moduli space converging to an $SL(m,\K)$ Higgs bundle gives just such a degeneration.

 \end{document}